%% file: int_quad_cuts.tex
\documentclass[12pt]{article}
\usepackage{fullpage,graphicx,psfrag,url}
\usepackage[small,bf]{caption}
\usepackage{amsmath,amssymb,amsthm,enumitem}
\usepackage{multirow}
\usepackage{tikz}
\usetikzlibrary{calc}
\input defs.tex

\newcounter{algorithmctr}[section]
\renewcommand{\thealgorithmctr}{\thesection.\arabic{algorithmctr}}
\newenvironment{algdesc}
   {\refstepcounter{algorithmctr}\begin{list}{}{
       \setlength{\rightmargin}{0\linewidth}
       \setlength{\leftmargin}{.05\linewidth}}
       \rmfamily\small
       \item[]{\setlength{\parskip}{0ex}\hrulefill\par
        \nopagebreak{\bfseries\textsf{Algorithm \thealgorithmctr~}}}}
   {{\setlength{\parskip}{-1ex}\nopagebreak\par\hrulefill} \end{list}}

\title{\Large{Concave Quadratic Cuts for Mixed-Integer Quadratic Problems}}
\author{Jaehyun Park \and Stephen Boyd}

\begin{document}
\maketitle

\begin{abstract}
The technique of semidefinite programming (SDP) relaxation can be used to
obtain a nontrivial bound on the optimal value of a nonconvex quadratically
constrained quadratic program (QCQP). We explore concave quadratic
inequalities that hold for any vector in the integer lattice $\integers^n$,
and show that adding these inequalities to a mixed-integer nonconvex QCQP can
improve the SDP-based bound on the optimal value. 
This scheme is tested using
several numerical problem instances of the
max-cut problem and the integer least squares problem.
\end{abstract}

\section{Introduction}

We consider mixed-integer indefinite quadratic optimization problems
of the form
\BEQ\label{problem-statement}
\begin{array}{ll}
\mbox{minimize} & f_0(x) = x^T P_0 x + q_0^T x + r_0 \\
\mbox{subject to} & f_i(x) = x^T P_i x + q_i^T x + r_i \le 0, \quad i = 1, \ldots, m \\
& x \in \mathcal{S},
\end{array}
\EEQ
with variable $x \in \reals^n$, where $\mathcal{S} = \{ x \,|\, x_1, \ldots,
x_p \in \integers \}$ is the \emph{mixed-integer set}, \ie, the set of 
real-valued vectors whose first 
$p$ components are integer-valued.
The problem data are $P_i \in \symm^n$, $q_i \in \reals^n$, and $r_i \in \reals$.
Here, $\symm^n$ denotes the set of $n \times n$ real-valued, symmetric,
possibly indefinite, matrices.
Quadratic equality constraints of the form $x^T F x + g^T x + h = 0$ can
also be handled by expressing them as two inequalities,
\[
x^T F x + g^T x + h \le 0, \quad -x^T F x - g^T x - h \le 0.
\]

The class of problems that can be written in the form
of~(\ref{problem-statement}) is very broad;
it includes other problem classes such as
mixed-integer linear programs (MILPs) and mixed-integer quadratic programs
(MIQPs).
Discrete constraints such as Boolean variables can be easily
encoded as well, which makes many NP-hard combinatorial optimization problems
special cases of~(\ref{problem-statement}). For example, $x_i^2 = 1$ encodes
the Boolean constraint that $x_i$ is either $+1$ or $-1$. The max-cut
problem is a well-known NP-hard problem with Boolean constraints only, which
can be formulated as the following:
\BEQ\label{maxcut-statement}
\begin{array}{ll}
\mbox{maximize} & -(1/4) x^T W x + (1/4) \ones^T W \ones
= (1/2) \sum_{i < j} W_{ij}(1 - x_i x_j) \\
\mbox{subject to} & x_i^2 = 1, \quad i=1, \ldots, n.
\end{array}
\EEQ
Here, $\ones$ represents a vector with all components equal to one. The matrix
$W \in \symm^n$ stores the weights of the edges in the graph:
$W_{ij}$ is the weight of the edge between node $i$ and $j$.
Other combinatorial problems in the form of~(\ref{problem-statement}) include
the maximum clique problem, graph bisection problem, and satisfiability
problem (SAT). In fact, any problem that is known to be an instance of
the quadratic assignment problem (QAP) fits into our
framework~\cite{koopmans1957assignment}.

Since~(\ref{problem-statement}) is an integer problem at its root, classical
number theoretic problems such as linear and quadratic Diophantine equations
are also special cases of mixed-integer
indefinite QCQP~\cite[\S 6]{nagell1951introduction}. The integer least
squares problem is another simple example of an integer quadratic problem:
\BEQ\label{ils-statement}
\begin{array}{ll}
\mbox{minimize} & \|Ax - b\|_2^2 \\
\mbox{subject to} & x \in \integers^n,
\end{array}
\EEQ
with variable $x$ and data $A \in \reals^{r \times n}$ and $b \in \reals^r$.
The integer least squares problem captures the essence of the phase ambiguity
estimation problem arising in the global positioning systems
(GPS)~\cite{hassibi1998integer}.

There are other interesting constraints that can be encoded in problems of the
form~(\ref{problem-statement}). For example, the
rank constraint $\Rank(X) \le k$ can be handled by introducing auxiliary matrix
variables $U$ and $V$ of appropriate dimensions, and adding a constraint
$X = UV$. Constraints involving the Euclidean distance between two points are
encoded naturally as well. The sphere packing problem~\cite{conway2013sphere}
and its variants are examples of problems involving (nonconvex) distance constraints.

Generic methods such as branch-and-bound~\cite{land1960automatic} or
branch-and-cut~\cite{padberg1991branch} can be used to
solve~(\ref{problem-statement}) globally, but they all have an exponential time
complexity. A more practical approach is to find an approximate solution,
or to compute lower and upper bounds on the optimal value.
The focus of this paper is on attaining a lower bound on the optimal
value of~(\ref{problem-statement}) by forming a \emph{semidefinite relaxation}
that is solvable in polynomial time. Semidefinite programming (SDP) is a
generalization of linear programming to symmetric positive semidefinite
matrices. The idea of semidefinite relaxation can be traced back
to~\cite{lovasz1979shannon}. Semidefinite relaxation is a powerful tool that
has been used not just within the domain of combinatorial
problems~\cite{luo2010semidefinite};
a notable example is sum-of-squares (SOS) optimization in control
theory~\cite{nesterov2000squared,parrilo2000structured}.
(For a thorough discussion of semidefinite programming and applications,
readers are directed to~\cite{vandenberghe1996semidefinite}.)
A well-known application of semidefinite relaxation is to
the max-cut problem; Goemans and Williamson constructed
a randomized algorithm for the
max-cut problem that attains a data-independent approximation factor of
$0.87856$, using the solution of a semidefinite
relaxation~\cite{goemans1995maxcut}.

Our main idea resembles that of cutting-plane method~\cite{kelley1960cutting},
and is closely related to hierarchies of linear and semidefinite programs
suggested by Lov{\'a}sz and Schrijver~\cite{lovasz1991cones}, Sherali and
Adams~\cite{sherali1990hierarchy}, Gomory~\cite{gomory1958outline},
Chv{\'a}tal~\cite{chvatal1973edmonds}, and
Lasserre~\cite{lasserre2001explicit}. However, to the best of our knowledge,
tightening SDP relaxation using true indefinite quadratic inequalities outside
the domain of Boolean problems is a novel approach.

\section{Semidefinite relaxation}\label{s-sdp}

A nontrivial lower bound on the optimal value of a (possibly nonconvex) QCQP
can be obtained by ``lifting'' it to a higher-dimensional space, and solving
the resulting problem. Consider a QCQP:
\BEQ\label{qcqp-statement}
\begin{array}{ll}
\mbox{minimize} & f_0(x) = x^T P_0 x + q_0^T x + r_0 \\
\mbox{subject to} & f_i(x) = x^T P_i x + q_i^T x + r_i \le 0, \quad i = 1, \ldots, m,
\end{array}
\EEQ
with variable $x \in \reals^n$. Let $f^\star$ denote its optimal value (which
can be $-\infty$). By introducing a new variable $X = x x^T$, we can
reformulate~(\ref{qcqp-statement}) as:
\[
\begin{array}{ll}
\mbox{minimize} & F_0(X, x) = \Tr(P_0 X) + q_0^T x + r_0 \\
\mbox{subject to} & F_i(X, x)
= \Tr(P_i X) + q_i^T x + r_i \le 0, \quad i=1, \ldots, m \\
& X = x x^T,
\end{array}
\]
with variables $X \in \symm^n$ and $x \in \reals^n$.

Then, we relax the nonconvex constraint $X = x x^T$ into a convex constraint
$X \succeq x x^T$ (where the operator $\succeq$ is with respect to the positive
semidefinite cone) and write it using a Schur complement to obtain a convex
relaxation:
\BEQ\label{sdp-relaxation}
\begin{array}{ll}
\mbox{minimize} & F_0(X, x) \\
\mbox{subject to} & F_i(X, x) \le 0, \quad i=1, \ldots, m \\
& \left[ \begin{array}{cc}
X & x \\ x^T & 1
\end{array} \right] \succeq 0.
\end{array}
\EEQ
This is now a convex problem in the ``lifted'' space, in fact an SDP
since the objective function is affine in $X$ and $x$, and its optimal
value is a lower bound on $f^\star$.

This SDP can be solved using an interior point method in polynomial time,
and in practice, the
number of iterations required is constant and insensitive to the problem size,
despite its worst-case complexity~\cite{vandenberghe1996semidefinite}.
A detailed analysis of the running time is beyond the scope of this paper.

\subsection{Lagrangian dual problem}

Here, we derive the Lagrangian dual problem~\cite[\S 5]{boyd2004convex} of
the SDP~(\ref{sdp-relaxation}). The Lagrangian of~(\ref{sdp-relaxation}) is
\[
L(X, x, \lambda, Y, y, \alpha) =
F_0 (X, x) + \sum_{i=1}^m \lambda_i F_i (X, x)
- \Tr \left[ \begin{array}{cc} Y & y \\ y^T & \alpha \end{array} \right]
\left[ \begin{array}{cc} X & x \\ x^T & 1 \end{array} \right].
\]
Let $g(\lambda, Y, y, \alpha)$ be the corresponding dual function, defined by
\[
g(\lambda, Y, y, \alpha) = \inf_{X, x} L(X, x, \lambda, Y, y, \alpha).
\]
The Lagrangian is affine in both $X$ and $x$, and thus minimizing it over $X$
and $x$ gives $-\infty$, unless the coefficients of $X$ and $x$ are both zero.
Therefore,
\[
g(\lambda, Y, y, \alpha) = \left\{ \begin{array}{ll}
r_0 + \sum_{i=1}^m \lambda_i r_i - \alpha
& Y = P_0 + \sum_{i=1}^m \lambda_i P_i, \quad y
= (1/2)(q_0 + \sum_{i=1}^m \lambda_i q_i) \\
-\infty & \mbox{otherwise.}
\end{array}\right.
\]
The dual problem is then
\[
\begin{array}{ll}
\mbox{maximize} & r_0 + \sum_{i=1}^m \lambda_i r_i - \alpha \\
\mbox{subject to} & Y = P_0 + \sum_{i=1}^m \lambda_i P_i \\
& y = (1/2)(q_0 + \sum_{i=1}^m \lambda_i q_i) \\
& \lambda_i \ge 0, \quad i=1, \ldots, m \\
& \left[ \begin{array}{cc}
Y & y \\ y^T & \alpha
\end{array} \right] \succeq 0,
\end{array}
\]
with variables $Y \in \symm^n$, $y \in \reals^n$, $\lambda \in
\reals^m$, $\alpha \in \reals$. Alternatively, by eliminating $Y$ and $y$,
we get an equivalent formulation
\BEQ\label{sdp-dual}
\begin{array}{ll}
\mbox{maximize} & r_0 + \sum_{i=1}^m \lambda_i r_i - \alpha \\
\mbox{subject to} & \lambda_i \ge 0, \quad i=1, \ldots, m \\
& \left[ \begin{array}{cc}
P_0 + \sum_{i=1}^m \lambda_i P_i
& (1/2)(q_0 + \sum_{i=1}^m \lambda_i q_i) \\
(1/2)(q_0 + \sum_{i=1}^m \lambda_i q_i)^T
& \alpha \end{array} \right] \succeq 0,
\end{array}
\EEQ
with variables $\lambda \in \reals^m$, $\alpha \in \reals$. Notice that the
dual problem~(\ref{sdp-dual}) is also a semidefinite program. Under mild
assumptions (\eg, strict feasibility of the primal problem), strong duality holds and
both~(\ref{sdp-relaxation}) and~(\ref{sdp-dual}) yield the same optimal value.

An advantage of considering this dual problem is that
unlike~(\ref{sdp-relaxation}), any feasible point
of~(\ref{sdp-dual}) yields a lower bound on
$f^\star$. This observation could be particularly useful if a dual feasible
solution with high objective value can be obtained relatively quickly.

\section{Concave quadratic cuts for mixed-integer vectors}\label{s-cuts}

\begin{figure}
  \centering
  \begin{tikzpicture}
    \coordinate (XMin) at (-4.5, 0);
    \coordinate (XMax) at (6.5, 0);
    \coordinate (YMin) at (0, -2.5);
    \coordinate (YMax) at (0, 4.5);
    \coordinate (APerp) at (-1, 1);
    \coordinate (E1) at (1, 0);
    \draw [thin, gray, -latex] (XMin)--(XMax);
    \draw [thin, gray, -latex] (YMin)--(YMax);

    \clip ($(XMin)+(YMin)$) rectangle ($(XMax)+(YMax)$);
    \foreach \x in {-4,-3,...,5}{
      \foreach \y in {-3,-2,...,4}{
        \node[draw,circle,inner sep=2pt,fill] at (2*\x, 2*\y) {};
      }
    }
    \filldraw[fill=gray, fill opacity=0.3, draw=black] ($10*(APerp)$)--($10*(APerp)+2*(E1)$)--($-10*(APerp)+2*(E1)$)--($-10*(APerp)$)--cycle;
  \end{tikzpicture}
  \caption{A concave quadratic inequality $(x_1+x_2)(x_1+x_2-1) \ge 0$. Notice
  that the shaded area, which does not satisfy the inequality, contains no
  lattice point.}
  \label{fig-int-quad-cuts}
\end{figure}
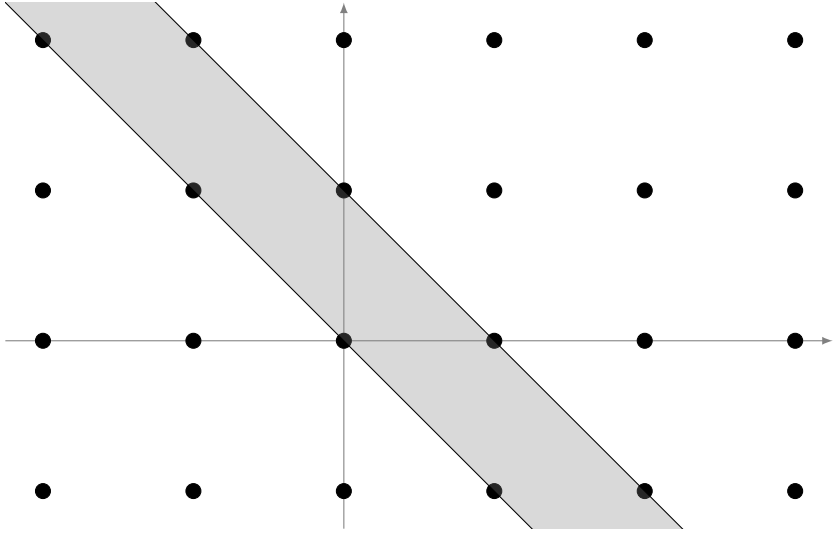

In this section, we present a simple method of relaxing the mixed-integrality
constraint $x \in \mathcal{S}$ into a set of concave quadratic inequalities.
Let $a \in \integers^n$ and $b \in \integers$ such that
$a_{p+1} = \cdots = a_n = 0$. The concave quadratic inequality
\BEQ\label{int-quad-cuts}
(a^T x - b)(a^T x - (b+1)) \ge 0,
\EEQ
or equivalently,
\[
-x^T (aa^T) x + (2b+1) a^T x - b(b+1) \le 0,
\]
holds if and only if $a^T x - b \le 0$ or $a^T x - b \ge 1$. In
particular,~(\ref{int-quad-cuts}) holds for \emph{every} vector $x \in
\mathcal{S}$, since then $a^T x - b$ is integer-valued, which (trivially)
satisfies $a^T x - b \le 0$ or $a^T x - b \ge 1$.
Figure~\ref{fig-int-quad-cuts} shows an
example of an inequality of the form~(\ref{int-quad-cuts}).
It follows from this observation that any number of such
inequalities can be added to~(\ref{problem-statement}) without changing the
optimal solution. That is, the following problem is equivalent to~(\ref
{problem-statement}):
\BEQ\label{equiv-statement}
\begin{array}{ll}
\mbox{minimize} & f_0(x) \\
\mbox{subject to} & f_i(x) \le 0, \quad i=1, \ldots, m \\
& - x^T (a_i a_i^T) x + (2b_i + 1) a_i^T x - b_i(b_i+1) \le 0, \quad i=1, \ldots, r \\
& x \in \mathcal{S},
\end{array}
\EEQ
with $a_i \in \integers^n$, $b_i \in \integers$, $(a_i)_{p+1} = \cdots =
(a_i)_n = 0$ for $i=1, \ldots, r$. Then, simply dropping the
mixed-integrality constraint from~(\ref{equiv-statement})
gives a nonconvex QCQP over $\reals^n$, and now the semidefinite relaxation
technique in~\S\ref{s-sdp} is readily applicable. The resulting SDP is:
\BEQ\label{relaxation-statement}
\begin{array}{ll}
\mbox{minimize} & F_0(X, x) \\
\mbox{subject to} & F_i(X, x) \le 0, \quad i=1, \ldots, m \\
& -\Tr(a_i a_i^T X) + (2b_i + 1) a_i^T x - b_i(b_i+1) \le 0, \quad i=1, \ldots, r \\
& \left[ \begin{array}{cc} X & x \\ x^T & 1 \end{array} \right] \succeq 0,
\end{array}
\EEQ
where each $F_i$ is defined as in~\S\ref{s-sdp}. The optimal value
$f^\mathrm{sdp}$ of~(\ref{relaxation-statement}) is a lower bound on
$f^\star$. Moreover, adding more true inequalities of the
form~(\ref{int-quad-cuts}) can only increase $f^\mathrm{sdp}$, hence
tightening the bound.

Inequalities of the form~(\ref{int-quad-cuts}) are satisfied by every point
in $\mathcal{S}$, but additionally, the mixed-integrality constraint
$x \in \mathcal{S}$ itself can be written as
a set of countably many such inequalities, \ie, the set
\[
\{x \,|\, (a^T x - b)(a^T x - (b+1))
\ge 0 \mbox{ for all } a \in \integers^n,
b \in \integers, a_{p+1} = \cdots = a_n = 0 \}
\]
is precisely $\mathcal{S}$. Thus, adding concave quadratic cuts can also be
interpreted as relaxing the mixed-integrality constraint, not by removing
it completely, but by replacing it with infinitely many concave quadratic
constraints, then dropping all but finitely many of them.

Although $\mathcal{S}$ can be written as a set of countably many concave
quadratic inequalities, when these inequalities are relaxed and rewritten in
the lifted space as in~(\ref{relaxation-statement}), the resulting set of
feasible points has no relationship with $\mathcal{S}$, in general.
Consequently, the solution or the optimal value
of~(\ref{relaxation-statement}) do not have any relationship
with $x^\star$ or $f^\star$. To demonstrate this, consider the following
problem in $\integers^2$:
\BEQ\label{concave-ex-problem}
\begin{array}{ll}
\mbox{minimize} & -\|x\|_2^2 \\
\mbox{subject to} & \|x\|_2^2 \le 1.2 \\
& x \in \integers^2,
\end{array}
\EEQ
which clearly has four optimal points $(\pm 1, 0)$ and $(0, \pm 1)$ with
objective value $f^\star = -1$. The SDP relaxation of this nonconvex problem
(without any additional concave quadratic cut) is
\[
\begin{array}{ll}
\mbox{minimize} & -\Tr X \\
\mbox{subject to} & \Tr X \le 1.2 \\
& \left[ \begin{array}{cc} X & x \\ x^T & 1 \end{array} \right] \succeq 0,
\end{array}
\]
with optimal value $f^\mathrm{sdp} = -1.2$.
Take $\hat{X} = 0.6I$ and $\hat{x} = 0$, which attain this objective value.
Then, for all $a \in \integers^n$ and $b \in \integers$,
\[
-\Tr(aa^T \hat{X}) + (2b + 1) a^T \hat{x} - b(b+1)
= -0.6\|a\|_2^2 - b(b+1) \le 0
\]
holds. The inequality follows from the fact that $b$ is integer-valued.
This shows that adding any number of inequalities in the form
of~(\ref{int-quad-cuts}) to~(\ref{concave-ex-problem}) will not
increase $f^\mathrm{sdp}$.

The example above shows that adding all possible concave quadratic cuts to an
integer problem and subsequently solving the relaxation, in general, does not
solve the original integer problem. In the special case of the integer least
squares problem, however, it is not clear whether this relaxation is tight or
not. When $n = 1$, \ie, when the problem reduces to minimizing a convex
quadratic function over the integers, it is easy to show that the relaxation
is indeed tight. On the other hand, even when $n = 2$, we have failed to either
prove the tightness of relaxation, or disprove it by producing any numerical
instance whose relaxation is not tight.

Finally, we discuss simple extensions of concave quadratic cuts.
The structure of the integer lattice, inherently, was not what made the
construction of~(\ref{int-quad-cuts}) possible. Rather, it was the
property that every feasible point of~(\ref{problem-statement})
satisfies \emph{exactly one} of the two affine inequalities
(namely $a^T x - b \le 0$ and $a^T x - (b+1) \ge 0$).
In other words, it does not matter whether the inequalities are
of the form $a^T x - b \le 0$ or not;
as long as there are two inequalities such that only
one of them holds for every feasible point, they can be multiplied together
to produce a true (and possibly nonconvex) inequality constraint. This
construction resembles that
of~\cite{shor1987quadratic,anstreicher2009semidefinite} in that affine
constraints are combined to quadratic constraints, and subsequently lifted
to a higher-dimensional space. However, the key difference is
that~(\ref{int-quad-cuts}) encodes \emph{exclusive disjunction}, \ie, the
two inequalities that we combine do not hold individually. There are other
types of concave quadratic inequalities that utilize the structure of the
integer lattice more. Take, for example, the constraint
\[
\|x - (1/2)\ones\|_2^2 \ge n/4.
\]
This constraint holds for every $x \in \integers^n$, but is not representable
as exclusive disjunction of any two affine inequalities. It is still a concave
quadratic cut, and thus can be used in our framework without any modification.
It is an open question as to whether these extensions result in a noticeable
improvement of our SDP bound.

\subsection{Choosing suitable cuts}\label{s-cut-choice}

Let $(\hat{X}, \hat{x})$ be a solution of~(\ref{relaxation-statement}),
and consider adding an additional inequality
\BEQ\label{additional-ineq}
-\Tr(aa^T X) + (2b + 1) a^T x - b(b+1) \le 0
\EEQ
to~(\ref{relaxation-statement}). Suppose that we want to choose integer-valued
$a$ and $b$ so that adding~(\ref{additional-ineq}) increases the SDP-based
lower bound $f^\mathrm{sdp}$. Then, we need $(\hat{X}, \hat{x})$ to violate
this inequality, \ie,
\BEQ\label{additional-ineq2}
-\Tr(aa^T \hat{X}) + (2b + 1) a^T \hat{x} - b(b+1) > 0.
\EEQ
Note that if $a$ is given, then choosing $b$ is easy, as we can maximize the
lefthand side with respect to $b$, and round it to the nearest integer value.
That is, we choose $b = \left\lfloor a^T \hat{x} \right\rfloor$.

Choosing a suitable vector $a$, unfortunately, is a difficult integer problem
itself, and we need to use a heuristic to find $a$. One very simple method is
to try all possible $a$ with at most $k$ nonzero elements that are either $+1$
or $-1$. For a given such vector $a$, it takes $O(k^2)$ time to check
whether~(\ref{additional-ineq2}) holds or not. Checking all possible such vectors
then takes $O(n^k 2^k k^2)$ time. In particular, if $k = 2$, then checking all
possible vectors can be done in $O(n^2)$ time.

A more involved heuristic uses an eigendecomposition of a certain matrix. To
derive the heuristic, we first bound the lefthand side
of~(\ref{additional-ineq2}) from below:
\BEAS
-\Tr(aa^T \hat{X}) + (2b + 1) a^T \hat{x} - b(b+1)
&\ge& \inf_{b \in \mathbf Z} -\Tr(aa^T \hat{X}) + (2b + 1) a^T \hat{x} - b(b+1)  \\
&\ge& \inf_{b \in \mathbf R} -\Tr(aa^T \hat{X}) + (2b + 1) a^T \hat{x} - b(b+1) \\
&=& -\Tr(aa^T \hat{X}) + 2 (a^T \hat{x})^2 - ((a^T \hat{x})^2 - 1/4) \\
&=& -\Tr(aa^T \hat{X}) + (a^T \hat{x})^2 + 1/4 \\
&=& - a^T (\hat{X} - \hat{x} \hat{x}^T) a + 1/4.
\EEAS
Our heuristic is to find $a \in \integers^n$ that maximizes the last line.
Recall that $M = \hat{X} - \hat{x} \hat{x}^T$ is positive semidefinite, and
thus $a = 0$ is a trivial maximizer of the expression. However, choosing $a =
0$ is not an option, since then~(\ref{additional-ineq2}) can never be
satisfied, regardless of $b$. Therefore, instead of choosing $a = 0$, we choose
an integer vector $a$ that is ``close'' to the eigenvector corresponding to
the smallest eigenvalue $v$ of $M$. At the same time, we want $\|a\|_2$ to be
small (but nonzero), because scaling $a$ by a factor of $t$ also scales $a^T M
a$ by a factor of $t^2$. Note that once $a$ is chosen this way, $b$ is set as
$\left\lfloor a^T \hat{x} \right\rfloor$ (which is the minimizer of the
lefthand side of~(\ref{additional-ineq2}) over $b \in \integers$), instead of
$a^T \hat{x} - 1/2$ (which is the minimizer of the lefthand side of~(\ref
{additional-ineq2}) over $b \in \reals$). Therefore,~(\ref{additional-ineq2})
may not hold even when $-a^T M a + 1/4 > 0$.
Conversely, if $-a^T M a + 1/4 \le 0$,
then~(\ref{additional-ineq2}) is guaranteed not to hold. In particular, if
$\lambda_{\min}$, the smallest eigenvalue of $M$, is larger than or equal to
$1/4$, then there exists no inequality of the form~(\ref{additional-ineq})
that can increase $f^\mathrm{sdp}$. It is not hard to justify this: \[ a^T M a
\ge \|a\|_2^2 \lambda_{\min} \ge \lambda_{\min}. \] There are a number of
reasonable ways to find a ``short'' integer vector $a$ that is (approximately)
aligned with a given vector $v$. For example, one can take an arbitrary
scaling factor $t > 0$ and round each entry of $tv$ to find $a$.
Alternatively, we can fix some small $k$, take $k$ indices $i_1, \ldots, i_k$
that correspond to the entries of $v$ with the largest magnitudes, and set
$a_{i_j} = \sign(v_{i_j})$ for $j=1, \ldots, k$, while leaving the other
entries of $a$ as zeros. In particular, when $k=1$, it means that $a$
will be of the form $a = \pm e_i$ for some $i$. Without loss of
generality, assume that $a = e_i$.
According to this heuristic, we choose $b =
\left\lfloor a^T \hat{x} \right\rfloor = \left\lfloor \hat{x}_i
\right\rfloor$, and thus the inequality~(\ref{additional-ineq}) we add in is
\[ - X_{ii} + (2 \left\lfloor \hat{x}_i \right\rfloor + 1) x_i -
\left\lfloor \hat{x}_i \right\rfloor (\left\lfloor \hat{x}_i
\right\rfloor + 1) \le 0. \] The corresponding concave quadratic inequality of
the form~(\ref{int-quad-cuts}) is \[ (x_i - \left\lfloor \hat{x}_i
\right\rfloor) (x_i - (\left\lfloor \hat{x}_i \right\rfloor + 1)) \ge 0. \]
Similarly, when $k=2$, the corresponding inequalities would be on $x_i \pm
x_j$ for some $i \ne j$.

Finally, we note that the heuristics described above can be applied in an
iterative manner, as in the well-known cutting-plane
method~\cite{kelley1960cutting}. That is, once some number of additional
cuts are introduced, we can solve the resulting SDP and find more cuts
using the new solution.

\subsection{Application to branch-and-cut}

In this section, we restrict ourselves to pure integer problems only, \ie,
$p=n$. Then, the concave quadratic cuts~(\ref{int-quad-cuts}) can be used in a
general branch-and-cut scheme to obtain the global solution. 
The branch-and-cut framework is depicted in Algorithm~\ref{alg-bnc}.

\begin{algdesc}\label{alg-bnc}
\emph{Branch-and-cut algorithm.}
\begin{tabbing}
{\bf given} an optimization problem $\mathcal{P}$ of
the form~(\ref{problem-statement}) with $p=n$. \\*[\smallskipamount]
1.\ \emph{Initialize.} Add $\mathcal{P}$ to $\mathcal{T}$,
the list of \emph{active problems}. Let $f^\star := \infty$.\\
{\bf while} $\mathcal{T}$ is nonempty \\
\qquad \= 2.\ \emph{Select an active problem.}
Remove $\mathcal{P'}$ from $\mathcal{T}$.\\
\> 3.\ Solve the SDP relaxation of $\mathcal{P'}$ to get its 
solution $(\hat{X}, \hat{x})$ with optimal value $f^\mathrm{sdp}$.\\
\> 4.\ {\bf if} $f^\mathrm{sdp} \ge f^\star$, go back to 2.\\
\> 5.\ {\bf if} $\hat{x} \in \integers^n$, set
$f^\star := f^\mathrm{sdp}$ and go back to 2.\\
\> 6.\ \emph{Cut.} If any $a \in \integers^n$ and
$b \in \integers$ satisfying~(\ref{additional-ineq2}) is found,
add~(\ref{additional-ineq}) to $\mathcal{P}'$ and go back to 3. \\
\> 7.\ Create two problem instances $\mathcal{P}_1$
and $\mathcal{P}_2$, both identical to $\mathcal{P}'$.\\
\> 8.\ Take some $c \in \integers^n$, $d \in \integers$,
then add constraint $c^T x \le d$ to $\mathcal{P}_1$,
and $c^T x \ge d+1$ to $\mathcal{P}_2$.\\
\> 9.\ \emph{Branch.} Add $\mathcal{P}_1$ and
$\mathcal{P}_2$ to $\mathcal{T}$.
\end{tabbing}
\end{algdesc}

There are a number of technical conditions that need to be met in order for
Algorithm~\ref{alg-bnc} to even terminate. We omit these details, as most
problems of practical interest meet these conditions, such as a bounded
domain, or the existence of the global solution (which would be implied by the
former condition). For example, any $0$-$1$ program satisfies these
requirements.

The two crucial steps that affect the overall performance of the
branch-and-cut algorithm are Steps 6 and 8. In~\S\ref{s-cut-choice}, we
have discussed a heuristic for finding a suitable cut that can be used in
Step 6. Step 8, which is called the \emph{branching step}, is another
important step in the algorithm. Algorithm~\ref{alg-bnc} may not even
terminate if the branching step is not implemented carefully;
for example, if the algorithm takes the same $c$ and $d$ at Step 7 every
time, then it produces redundant branches, making no progress as a result. A
commonly used branching strategy is to branch on a particular variable $x_i$.
That is, for some index $i$ such that $\hat{x}_i$ is not integer-valued, we
add $x_i \le \left\lfloor \hat{x}_i \right\rfloor$ and
$x_i \ge \left\lfloor \hat{x}_i \right\rfloor + 1$ as the branching
inequalities. This would correspond to choosing $c = e_i$ and
$d = \left\lfloor c^T \hat{x} \right\rfloor$ in Step 8 of
Algorithm~\ref{alg-bnc}.

While branching on a single variable is an intuitive and simple strategy that
is commonly used, we generalize this idea and find a ``natural'' branching
inequality that can be easily obtained from the solution of the SDP
relaxation~(\ref{sdp-relaxation}), with almost no additional computation. The
main idea comes from sensitivity analysis of convex problems;
the optimal dual variables of~(\ref{sdp-relaxation}) gives information about
the sensitivity of $f^\mathrm{sdp}$, with respect to perturbations of the
corresponding constraints. Roughly speaking, if the magnitude of
$\lambda_i$---the optimal dual variable corresponding to the constraint
$F_i(X, x) \le 0$ of~(\ref{sdp-relaxation})---is big, then the constraint
is ``tight,'' and further tightening the constraint would lead to a
big increase in $f^\mathrm{sdp}$. To be precise, let $f^\mathrm{sdp}(u)$
denote the optimal value of problem~(\ref{sdp-relaxation}), when the
inequality constraint $F_i(X, x) \le 0$ is replaced with $F_i(X, x) \le u$.
Let $\lambda_i$ be the optimal dual variable corresponding to the constraint
of the unperturbed problem~(\ref{relaxation-statement}). Then, for all $u$, we
have
\[
f^\mathrm{sdp}(u) \ge f^\mathrm{sdp}(0) - \lambda_i u.
\]
In other words, tightening an inequality (\ie, $u < 0$) by $|u|$ increases the
SDP bound by at least $\lambda_i |u|$. According to this interpretation,
natural branching inequalities come from the concave quadratic cut (in the
lifted space) with the largest dual variable. That is, if the constraint
\BEQ\label{lifted-ineq}
-\Tr(a a^T X) + (2b + 1) a^T x - b(b+1) \le 0
\EEQ
has the largest value of the dual variable (\ie, is the tightest), then we add
$a^T x \le b$ and $a^T x \ge b + 1$ as branching inequalities. The intuition
behind this choice is that these inequalities are tighter versions
of~(\ref{lifted-ineq});
recall that~(\ref{lifted-ineq}) is a relaxation of~(\ref{int-quad-cuts}),
which is satisfied exactly when $a^Tx \le b$ or $a^T x \ge b+1$.
Therefore, we expect $f^\mathrm{sdp}$ to go up by adding any of these
inequalities. After a branching inequality is added,~(\ref{lifted-ineq}) can
be removed, as it is implied by the newly added branching inequality.
When there is no concave quadratic
cut with strictly positive dual variable, then we may branch on a single
variable.

We note that scaling a constraint by a factor of $t > 0$ scales the
corresponding optimal dual variable by a factor of $1/t$. Therefore, to
correctly compare the dual variables and add branching inequalities, we have
to apply a proper scaling to each constraint of the
form~(\ref{additional-ineq}).
It is difficult to determine the scaling factor directly from perturbation
and sensitivity analysis, as the branching inequalities we add, despite
being tighter than~(\ref{lifted-ineq}), are not direct perturbations
of it. However, experiments suggest that~(\ref{additional-ineq}) is already
properly scaled, and thus directly comparing the optimal dual variables gives
good branching inequalities.

\section{Examples}

In this section, we consider numerical instances of the integer least squares
problem and max-cut problem to show the effectiveness of the quadratic concave
cuts in terms of the SDP-based lower bound. We emphasize that these problems
were chosen because they are simple to describe, and showed qualitatively
different results. They are by no means representative problems of the entire
class of problems that can be written as~(\ref{problem-statement}).

\subsection{Computational details}
The SDP~(\ref{sdp-relaxation}) was solved using CVX~\cite{cvx,gb08} with the
MOSEK 7.1 solver~\cite{mosek}, on a 3.40 GHz Intel Xeon machine. In order to
obtain the solution in a reasonable amount of time, we only considered 
small-sized problems of $n \sim 100$.

\subsection{Integer least squares}

\paragraph{Problem formulation.}
We consider the following problem formulation that is equivalent
to~(\ref{ils-statement}):
\BEQ\label{ils-reformulation}
\begin{array}{ll}
\mbox{minimize} & \|A(x - x^\mathrm{cts})\|_2^2 \\
\mbox{subject to} & x \in \integers^n.
\end{array}
\EEQ
Here, $x^\mathrm{cts} \in \reals^n$ is a given point at which the objective
value becomes zero, which is a simple lower bound on the optimal
value $f^\star$.

\paragraph{Problem instances.}
We use random problem instances of the integer least squares
problem~(\ref{ils-reformulation}), generated in the same way as
in~\cite{park2015semidefinite}:
entries of $A \in \reals^{r \times n}$ are
sampled independently from $\mathcal{N}(0, 1)$, with the number of rows $r$
set as $r = 2n$. The point $x^\mathrm{cts}$ was drawn from the uniform
distribution on the box $[0, 1]^n$.

\paragraph{Method.}
We compare three lower bounds and one upper bound on the optimal value
$f^\star$. The first lower bound is the simple lower bound of the continuous
relaxation: $f^\mathrm{cts} = 0$. The second lower bound, which we denote by
$f^\mathrm{sdp}_1$, is the SDP-based lower bound explored
in~\cite{park2015semidefinite}. This was obtained by relaxing the integer
constraint to a set of $n$ quadratic concave inequalities $x_i (x_i - 1) \ge
0$ for all $i$, followed by solving the SDP relaxation. The third
lower bound, $f^\mathrm{sdp}_2$, was obtained by generalizing this approach,
as described in~\S\ref{s-cuts}; in addition to the inequalities
$x_i (x_i - 1) \ge 0$, we considered $O(n^2)$ additional cuts
with $\|a\|_2 = \sqrt{2}$, and added only those
satisfied~(\ref{additional-ineq2}). Adding more inequalities (with
$\|a\|_2 \ge \sqrt{3}$) gave little improvement, and thus they were not
considered for our experiments. The upper
bound, $\hat{f}$, was found by running the randomized algorithm
constructed from the solution of the SDP
(see~\cite{goemans1995maxcut,park2015semidefinite}).
In~\cite{park2015semidefinite}, it was shown empirically that this upper bound
is very close to the optimal value for problems of small enough size ($n \le
60$) that the optimal value was obtainable.

\paragraph{Results.}
For each problem size, we generated $100$ problem instances and collected the
lower and upper bounds from the SDP relaxation. In Table~\ref{t-ils-runtime},
we show the average running time of the two SDPs used to obtain
$f^\mathrm{sdp}_1$ and $f^\mathrm{sdp}_2$, respectively, and the average
number of additional cuts added to the second SDP. The trade-off between the
number of cuts (which was roughly $n^2/4$) and the running time of the SDP is
clear from the table.

\begin{table}
\begin{center}
\begin{tabular}{|c||c|c|c|c|}
\hline 
$n$ & SDP1 & SDP2 & \# of cuts \\ \hline\hline
40   & $0.378$ & $1.984$ & $516.2$ \\ \hline
50   & $0.338$ & $3.542$ & $749.3$ \\ \hline
60   & $0.374$ & $6.584$ & $1033$ \\ \hline
70   & $0.443$ & $11.82$ & $1380$ \\ \hline
80   & $0.543$ & $19.12$ & $1732$ \\ \hline
100  & $0.830$ & $47.51$ & $2579$ \\ \hline
\end{tabular}
\end{center}
\caption{Running time of SDPs by number of variables, and the average number of additional cuts added to the second SDP.}
\label{t-ils-runtime}
\end{table}

In Table~\ref{t-ils-bound-comp}, we compare the two lower bounds
$f^\mathrm{sdp}_1$ and $f^\mathrm{sdp}_2$, along with an upper bound
$\hat{f}$. The simple lower bound $f^\mathrm{cts} = 0$ was
omitted from the table. Note that $f^\mathrm{sdp}_2 \ge f^\mathrm{sdp}_1$
holds not only on average, but for every problem instance, because the second
SDP is more constrained than the first.
The ratio between the two optimality gaps, namely
\[
\alpha = \frac{\hat{f} - f^\mathrm{sdp}_2} {\hat{f} - f^\mathrm{sdp}_1},
\]
is also shown in the same table. We obtained a significant
reduction for all problem sizes ($\alpha = 0.61$ for $n=100$), though we
expect $\alpha$ to be higher for larger problems. Note, however, that
$\hat{f}$ is not the optimal value, and thus the true reduction in the
optimality gap is larger than what is reported. This is a notable
improvement, because~\cite{park2015semidefinite}, to the best of our
knowledge, is a state of the art method for obtaining a lower bound
on the integer least squares problem (in polynomial time).

\begin{table}
\begin{center}
\begin{tabular}{|c||c|c|c|c|}
\hline 
$n$ & $f^\mathrm{sdp}_1$ & $f^\mathrm{sdp}_2$ & $\hat{f}$ & $\alpha$ \\ \hline\hline
40   & $88.21$ & $131.9$ & $165.6$ & $0.43$ \\ \hline
50   & $135.0$ & $198.8$ & $259.4$ & $0.48$ \\ \hline
60   & $186.0$ & $272.2$ & $369.2$ & $0.53$ \\ \hline
70   & $245.1$ & $356.7$ & $493.4$ & $0.55$ \\ \hline
80   & $310.3$ & $451.1$ & $646.0$ & $0.58$ \\ \hline
100  & $469.9$ & $675.7$ & $1001$  & $0.61$ \\ \hline
\end{tabular}
\end{center}
\caption{Average lower and upper bounds by number of variables, along with the average ratio between the optimality gaps.}
\label{t-ils-bound-comp}
\end{table}

\subsection{Max-cut problem}

\paragraph{Problem formulation.}
We reformulate~(\ref{maxcut-statement}) in terms of $0$-$1$ variables $z_i =
(1/2)(x_i + 1)$, so that the cuts introduced in~\S\ref{s-cuts} are tighter:
\BEQ\label{maxcut-reformulation}
\begin{array}{ll}
\mbox{maximize} & (W \ones)^T z - z^T W z \\
\mbox{subject to} & z_i(z_i - 1) = 0, \quad i=1, \ldots, n.
\end{array}
\EEQ
Note that~(\ref{maxcut-reformulation}) is a maximization problem, and hence we
get an upper bound on the optimal value by solving its relaxation.

\paragraph{Problem instances.}
We used the set of $10$ small-sized problems ($n=125$) with $\pm 1$ edge
weights, which was used in~\cite{festa2002randomized}. Note that due to the
integral edge weights, the optimal value also is integer-valued. It follows
that if $f'$ is an upper bound on the optimal value $f^\star$, then
$\left\lfloor f' \right\rfloor$ is also an upper bound on $f^\star$. In
particular, if we have some feasible point $z$ and some relaxation of the
max-cut problem has an optimal value $f^\mathrm{sdp}$ such that $f_0(z)
\le f^\mathrm{sdp} < f_0(z) + 1$, then the optimal solution to the SDP
provides a \emph{certificate of optimality} of the point $z$, \ie, $f_0(z) =
f^\star$.

\paragraph{Method.}
In this application, we compare two upper bounds, and the best known lower
bound to each of the 10 problem instances. The first upper bound
$f^\mathrm{gw}$ is the classical SDP bound explored
in~\cite{goemans1995maxcut}. This bound was obtained by
relaxing~(\ref{maxcut-reformulation}) according to~\S2,
without adding any concave quadratic cut.
The second lower bound, which we denote by $f^\mathrm{sdp}$, was found by
adding concave quadratic cuts satisfying~(\ref{additional-ineq2}), just
as in the integer least squares problem. However, for the max-cut problem,
we found effectively no improvement in the SDP
bound by adding inequalities with $\|a\|_2 = \sqrt{2}$, even when all such
inequalities were added. Instead, we generated $10{,}000$ random vectors $a$
with $\|a\|_2 = \sqrt{3}$ (\ie, having exactly three nonzero entries that are
$\pm 1$), and added the corresponding cuts to the SDP relaxation only
when~(\ref{additional-ineq2}) was satisfied.

\paragraph{Results.}
In Table~\ref{t-maxcut-runtime}, we show the solve time of the two SDPs used
to obtain $f^\mathrm{gw}$ and $f^\mathrm{sdp}$, respectively, and the number
of additional cuts added to the second SDP. As seen in the previous
application, the trade-off between the number of cuts and the running time is
clear.

\begin{table}
\begin{center}
\begin{tabular}{|c||c|c|c|}
\hline 
Instance & GW & SDP & \# of cuts \\ \hline\hline
G54100   & 1.92 & 17.9 & 1072 \\ \hline
G54200   & 1.07 & 16.4 &  997 \\ \hline
G54300   & 1.07 & 13.1 &  875 \\ \hline
G54400   & 0.86 & 18.1 & 1117 \\ \hline
G54500   & 0.86 & 16.4 &  987 \\ \hline
G54600   & 1.17 & 14.1 &  915 \\ \hline
G54700   & 0.82 & 20.1 & 1115 \\ \hline
G54800   & 1.07 & 15.8 &  969 \\ \hline
G54900   & 0.91 & 23.0 & 1325 \\ \hline
G541000  & 0.89 & 22.9 & 1343 \\ \hline
\end{tabular}
\end{center}
\caption{Running time of two SDP relaxations, and the number of additional cuts added to the second SDP.}
\label{t-maxcut-runtime}
\end{table}

In Table~\ref{t-maxcut-bound-comp}, we compare the two upper bounds
$f^\mathrm{gw}$ and $f^\mathrm{sdp}$, along with the best known lower bound
$\hat{f}$. The upper bounds were rounded down to the nearest integer, using
the observation made above. The ratio between the optimality gaps
\[
\alpha = \frac{f^\mathrm{sdp} - \hat{f}} {f^\mathrm{gw} - \hat{f}},
\]
is also shown in the same table. The value of $\alpha$ ranged from $0.73$ to $0.88$,
which is larger than what we obtained in the case of integer least squares.

\begin{table}
\begin{center}
\begin{tabular}{|c||c|c|c|c|}
\hline 
Instance & $f^\mathrm{gw}$ & $f^\mathrm{sdp}$ & $\hat{f}$ & $\alpha$ \\ \hline\hline
G54100   & 126 & 123 & 110 & 0.81 \\ \hline
G54200   & 128 & 125 & 112 & 0.81 \\ \hline
G54300   & 123 & 121 & 106 & 0.88 \\ \hline
G54400   & 128 & 125 & 114 & 0.79 \\ \hline
G54500   & 127 & 123 & 112 & 0.73 \\ \hline
G54600   & 126 & 124 & 110 & 0.88 \\ \hline
G54700   & 126 & 124 & 112 & 0.86 \\ \hline
G54800   & 125 & 122 & 108 & 0.82 \\ \hline
G54900   & 126 & 123 & 110 & 0.81 \\ \hline
G541000  & 127 & 124 & 112 & 0.80 \\ \hline
\end{tabular}
\end{center}
\caption{Best known lower bound $\hat{f}$, Goemans-Williamson SDP bound $f^\mathrm{gw}$, our SDP-based upper bound $f^\mathrm{sdp}$, and the reduction in optimality gap for each problem instance.}
\label{t-maxcut-bound-comp}
\end{table}

It should be noted that the objective of this experiment is not to compete
with state of the art methods for solving the max-cut problem such
as~\cite{marti2009advanced} (which are shown to find near-optimal solutions
very efficiently), but to demonstrate that our
method generates a nontrivial upper bound for Boolean problems that is
better than the Goemans-Williamson bound. Whether this method has any
direct relationship with the (3rd) Lasserre hierarchy is an outstanding
question.

\nocite{*}\bibliography{refs}
\end{document}

%% file: defs.tex
\newcommand{\BEAS}{\begin{eqnarray*}}
\newcommand{\EEAS}{\end{eqnarray*}}
\newcommand{\BEQ}{\begin{equation}}
\newcommand{\EEQ}{\end{equation}}
\newcommand{\BIT}{\begin{itemize}}
\newcommand{\EIT}{\end{itemize}}

\newcommand{\eg}{{\it e.g.}}
\newcommand{\ie}{{\it i.e.}}

\newcommand{\ones}{\mathbf 1}
\newcommand{\reals}{{\mbox{\bf R}}}
\newcommand{\integers}{{\mbox{\bf Z}}}
\newcommand{\symm}{{\mbox{\bf S}}}  

\newcommand{\Rank}{\mathop{\bf Rank}}
\newcommand{\Tr}{\mathop{\bf Tr}}

\newcommand{\sign}{\mathop{\bf sign}}

\newcounter{oursection}